\magnification=\magstep1
\input amstex
\documentstyle {amsppt}
\NoBlackBoxes
\NoRunningHeads
\pageheight {8.5 truein}
\pagewidth {6.5 truein}

\topmatter
\title
On a class of algebras associated to directed graphs
\endtitle
\author
Israel Gelfand, Vladimir Retakh \footnote {partially supported by NSA \hfill \hfill}, Shirlei Serconek \footnote  {partially supported by CNPq/PADCT \hfill \hfill} 
 and Robert Lee Wilson
\endauthor
\address
\newline
I.G., V.R., R.W.: Department of Mathematics, Rutgers University,
Piscataway,
NJ 08854-8019
\newline
S.S:
IME-UFG
CX Postal 131
Goiania - GO
CEP 74001-970 Brazil
\endaddress
\email
\newline igelfand\@math.rutgers.edu, vretakh\@math.rutgers.edu, 
\newline serconek\@math.rutgers.edu, rwilson\@math.rutgers.edu
\endemail

\endtopmatter
\document

\centerline {\bf Abstract}
\smallskip To any directed graph we associate an algebra with
edges of the graph as generators and with relations
defined by all pairs of directed paths with the same
origin and terminus. Such algebras are related to factorizations
of polynomials over noncommutative algebras.
We also construct a basis for our algebras
asssociated to layered graphs.
\bigskip

\head 0. Introduction\endhead

Factorizations of noncommutative polynomials play an important
role in many areas of mathematics  such as operator theory,
integrable systems, and Yang-Baxter equations (see, for example,
\cite {O, V}). In this paper we use directed graphs and
algebras associated with those graphs as a natural framework for
studying such factorizations.

Let $R$ be an associative ring with unit and $P(\tau )\in R[\tau ]$
a polynomial over $R$, where $\tau$ is a central variable.
Assume that $P(\tau )$ is {\it a monic polynomial} (i.e. the leading coefficient
in $P(\tau )$ equals $1$). When $R$ is a (commutative) field there exists at most one
(up to rearrangement of the factors)  factorization of a monic polynomial $P(\tau )$ of degree $n$
into a product of linear polynomials:
$$P(\tau )=(\tau -y_n) (\tau -y_{n-1})\dots (\tau -y_1) \tag 0.1$$
When the ring $R$ is not commutative, there may exist many factorizations of
type (0.1).

In \cite {GGRSW} the elements $y_1, y_2,\dots , y_n$ in formula (0.1) where called
the {\it pseudoroots} of the polynomial $P(\tau)$. The element $y_1$  is {\it a right root} of $P(\tau)$
and the element $y_n$ is {\it a left root} of $P(\tau )$.

Let  $X=\{x_1, x_2, \dots , x_n\}\subseteq R$ be {\it a generic
set} of right roots of $P(\tau )$ (meaning that if $k \ge 2$ and 
$\{x_{i_1}, x_{i_2},\dots , x_{i_k}\} \subseteq X$, the corresponding
Vandermonde matrix is invertible). It was shown in
\cite {GR1} (see also \cite {GR2, GGRW}) that the right roots $x_1,
x_2, \dots , x_n$ define a set of pseudoroots $x_{A, i}$ of
$P(\tau)$. Here $A\subseteq \{1,2,\dots , n\}$, $i\in \{1,2,\dots ,
n\}\setminus A$, and $x_{\emptyset , j}=x_j$ for all $j$.

According to \cite {GR1} (see also \cite {GR2, GGRW}) for any ordering 
$i_1, 
i_2,\dots , i_n$ of
$1,2,\dots , n$ the pseudoroots $y_k=x_{\{i_1,\dots i_{k-1}\}, i_k}$,
$k=1,2,\dots , n$ define a factorization (0.1) .

The pseudoroots $x_{A,i}$ satisfy the following identities for any $i,j\notin A$
$$x_{A\cup \{i\}, j}+x_{A, i}= x_{A\cup \{j\}, i}+x_{A, j} \tag 0.2a $$
$$x_{A\cup \{i\}, j}x_{A, i}= x_{A\cup \{j\}, i}x_{A, j}. \tag 0.2b $$

The paper \cite {GRW}  introduces and studies the algebra
$Q_n$, called the universal algebra of pseudoroots, generated
by elements $x_{A, i}$ satisfying the identities (0.2).

The study of $Q_n$ arose from the
theory of quasideterminants  and, specifically, from the
noncommutative version of the Vi\`{e}te Theorem and its relations to
the theory of noncommutative symmetric functions (see, for
example, \cite {GGRW}). In particular, the algebra $Q_n$  is
quadratic, Koszul, and its dual algebra $Q_n^!$ has  finite
dimension (see \cite {GGRSW, SW}).

A natural description of the algebra $Q_n$ (and, therefore,
factorizations of noncommutative polynomials) can be given by
using directed graphs. Let $\Gamma _n$ be the Hasse graph
corresponding to the lattice of subsets of  the set $\{1,2,\dots
, n\}$. Vertices of $\Gamma _n$ are subsets $A\subseteq
\{1,2,\dots , n\}$. Edges of $\Gamma _n$ are defined by pairs
$(A, i)$ where $i\in \{1,2,\dots , n\}\setminus A$. Such edges
go from the vertex $A\cup \{i\}$ to the vertex $A$.

The graph  $\Gamma _n$ possesses several properties. It is a hypercube with
$2^n$ vertices. It is {\it a layered graph}: each vertex $A$ has {\it a level}
$|A|=\text {card} \ A$. There is only one vertex $\{1,2,\dots , n\}$ of 
level
$n$ and only one vertex $* $ of level $0$ and any edge $a$ goes
from a vertex of level $r$ to a vertex of level $r-1$.

One can describe the relations in the algebra $Q_n$  by using geometric properties of
the directed graph $\Gamma_n$. Fix a field $F$. Let $T$ be the free associative algebra over
$F$ generated by the edges of $\Gamma _n$.  Then the algebra $Q_n$ is a quotient algebra
of $T$ modulo relations (0.2).

To every directed path $\pi =(a_1, a_2, \dots , a_m)$ in $\Gamma _n$ (where $a_1,...,a_m$ 
are edges of $\Gamma_n$) there is a corresponding polynomial
$P_{\pi }(\tau )\in T[\tau ]$, $P_{\pi }(\tau ) =(\tau -a_1)(\tau -a_2)\dots (\tau -a_m)$.
Denote the image of $P_{\pi}(\tau)$ in $Q_n[\tau ]$ by $\tilde P_{\pi }(\tau )$.
If  two paths $\pi _1$ and $\pi _2$ both
go from a vertex $v$ to a vertex $v'$ then
$$\tilde P_{\pi _1}(\tau )=\tilde P_{\pi _2}(\tau ) \tag 0.3$$
in $Q_n[\tau ]$.

The identity (0.3) defines relations in $Q_n$. Note that the defining relations (0.2) arise from those 
pairs of paths forming a diamond, i. e. a path $\pi _1$ consisting of edges
going from $A\cup \{i,j\}$ to $A\cup \{i\}$ and from   $A\cup \{i\}$  to $A$,
and  a  path $\pi _2$ consisting of edges
going from $A\cup \{i,j\}$ to $A\cup \{j\}$ and from   $A\cup \{j\}$  to $A$.

In this paper we introduce and study a class of algebras associated to certain directed graphs
as a natural generalization of the geometric definition
of the  algebra $Q_n$. It is more convenient for us to use formulas similar to
(0.3) by introducing a new variable $t=\tau ^{-1}$.

Let $\Gamma =(V, E)$ be a directed graph with vertices $V$ and
edges $E$. Assume that $\Gamma$ is layered and that the maximum level is $n$.
Fix a field $F$ and denote by $T(E)$ the free
associative algebra over  $F$ generated by all edges. To any path
$\pi =(e_1, e_2,\dots , e_k)$ in $\Gamma $ corresponds a
polynomial
$$P_{\pi }(t)=(1-te_1)(1-te_2)\dots (1-te_k)\in T(E)[t]/(t^{n+1}).$$
 We define the algebra $A(\Gamma )$ to be the  quotient algebra of $T(E)$ modulo relations defined by
the equalities implied by
$$P_{\pi _1}(t)=P_{\pi _2}(t ) \tag 0.4 $$
 where the paths $\pi _1$ and $\pi _2$ both go
from a vertex $v$ to a vertex $v'$.

The algebra $A(\Gamma_n) $  coincides with the algebra $Q_n$.

Our main examples of directed graphs are the Hasse graphs associated with lattices of
subsets and subspaces, abstract polytopes, complexes, and partitions.

We have already mentioned that the algebra $Q_n$, its subalgebras,
quotient algebras, and dual algebras  possess many interesting properties (see \cite {GRW, GGR, GGRSW,
SW, Pi}). We believe that the same is true for  algebras associated to other directed graphs. In particular,
we are planning to relate the structure of the algebra $A(\Gamma )$ with the geometry of the graph.

As a first step in these directions, we construct in this
paper a linear basis for the  algebra $A(\Gamma )$ for a  large class
of directed graphs $\Gamma $. This basis has a ``geometric"
nature. In particular,
our result  simplifies  the construction of the basis for $Q_n$
given in  \cite {GRW}.

We  now briefly describe this basis.  For each vertex $v \in V$ with ${\text {level}} \ |v| > 0$, we choose (arbitrarily) an edge $e_v$ beginning at the vertex $v$.  This defines a unique path $\pi_v$ from $v$ to $*$.  In $A(\Gamma)$ write
$\tilde{e}$ for the image of $e \in E$ and
$$P_{\pi_v}(t) = \sum_{k=0}^{|v|} \tilde {e}(v,k)t^k.$$
We say that two pairs $(v,k), (w,m) \in V \times {\bold Z}_{\ge 0}$ can be composed (and write
$(v,k) \models (w,m)$) if there is a path of length $k - m$ from $v$ to $w$.  
Let ${\bold B}(\Gamma)$   be the set of all sequences
$$ {\bold b} = ((b_1,m_1),(b_2,m_2),...,(b_k,m_k))$$  where 
$k \geq 0$, $b_1, b_2, ..., b_k \in V$, $1 \leq m_i \leq |b_i|$ for $1 \leq i \leq k$, and
$(b_i,m_i) \not\models (b_{i+1},m_{i+1})$ for $1 \le i < k$.

For ${\bold b} = ((b_1,m_1),...,(b_k,m_k)) \in {\bold B}(\Gamma)$ set
$$\tilde {e}({\bold b}) = \tilde {e}(b_1,m_1)... \tilde {e}(b_k,m_k).$$

\proclaim{Theorem 4.3} Let $\Gamma = (V,E)$  be a layered graph, $V = \cup_{i=0}^n V_i, $ and $V_0 = \{*\}$ 
where $*$ is the unique minimal vertex of $\Gamma$. Then $\{ \tilde{e}({\bold b}) \ | \  {\bold b} \in {\bold B}(\Gamma) \}$ is a basis 
for $A(\Gamma)$.
\endproclaim

An equivalent formulation of this theorem may be obtained by replacing each ``generating function coefficient" 
$\tilde{e}(v,k)$ by the ``monomial" $ \tilde{e}_{v^{(0)}}\tilde{e}_{v^{(1)}}...\tilde{e}_{v^{(k-1)}}$
where 
$v^{(0)}$, $v^{(1)}$, $\dots $, $v^{(l)} = *$ are the vertices of the path 
$\pi_v$.  (This monomial is the leading term of 
$\tilde{e}(v,k)$ in an appropriate filtration.)  
This ``monomial" formulation of our basis theorem, when specialized to the 
graph $\Gamma_n$, gives the basis theorem of \cite {GRW} for $Q_n$.
\bigskip
\head \bf 1. The directed graph $\Gamma = (V, E)$ \endhead
\bigskip
Let  $\Gamma = (V, E)$ be a {\bf directed graph}. That is, $V$ is a set (of vertices), 
$E$ is a set (of edges), and ${\bold t}: E \rightarrow V$ and ${\bold h}: E \rightarrow V$ are functions. (${\bold t}(e)$ is the {\it tail} 
of $e$ and ${\bold h}(e)$ is the {\it head} of $e$.)       
\smallskip
We say that $\Gamma$ is {\bf layered} if $V = \cup _{i=0}^n V_i$, $E = \cup_{i=1}^n E_i$, 
 ${\bold t}: E_i \rightarrow V_i$, \ ${\bold h}: E_i \rightarrow V_{i-1}$.
\smallskip  
We will assume throughout the remainder of the paper that $\Gamma = (V, E)$ is a layered graph with 
$V = \cup_{i=0}^n V_i$, and $V_0 = \{*\}$ where $*$ is the unique minimal vertex of $\Gamma$ (i.e., for every 
$v \in V, v \ne *$, there exists $e \in E$ with ${\bold t}(e) = v$).

For each $v \in \cup_{i=1}^n V_i$ we will fix, arbitrarily, 
some $e_v \in E$, with ${\bold t}(e) = v$.
If $v \in V_i$ we write $|v| = i$ and say that $v$ has {\it level} $i$. Similarly, 
if $e \in E_i$ we write $|e| = i$ and say that $e$ has {\it level} $i$.

If $v, \ w \in V$, a {\bf path} from $v$ to $w$ is a sequence of edges 
$\pi = \{ e_1, e_2, ...,e_k \}$ with ${\bold t}(e_{1}) = v$, ${\bold h}(e_k) = w$ and ${\bold t}(e_{i+1}) = {\bold h}(e_i)$ for 
$1 \leq i < k$.  We write $v = {\bold t}(\pi)$, $w = {\bold h}(\pi)$.  We also write $v > w$ if there  is a path from 
$v$ to $w$.

Let $l(\pi)$, the {\it length} of $\pi$, denote $k$, and let $|\pi|$, the {\it level} of $\pi$, denote 
$|e_1| + ... + |e_k|$.
 
If $\pi_1 = \{ e_1, ...,e_k \}$,$\pi_2 = \{ f_1, ...,f_l \}$ are paths with ${\bold h}(\pi_1) = {\bold t}(\pi_2)$ 
then $\{ e_1, ...,e_k, f_1, ..., f_l \}$ is a path; we denote it by $\pi_1\pi_2$.
      
For $v \in V$, write $v^{(0)} = v$ and define $v^{(i+1)} = {\bold h}(e_{v^{(i)}})$ for $0 \leq i <  |v|$. 
Then $v^{(|v|)} = *$ and $\pi_v = \{ e_{v^{(0)}}, ...,e_{v^{(|v|-1)}} \}$  is a path from $v$ to $*$. 

\bigskip
\head \bf 2. The filtered algebra $T(E)$
\endhead
\bigskip
Let $T(E)$ denote the free associative algebra  on $E$ over a field $F$. Define
$$T(E)_i = span \{ e_1...e_r \ | \ r \geq 0, |e_1| + ... + |e_r| \leq i \}.$$

If $a \in T(E)_i$, $a \not\in T(E)_{i-1}$, write $|a| = i$.
For a path $\pi = \{ e_1, e_2, ...,e_k \}$ define $$P_{\pi}(t) = (1-te_1)...(1-te_k) \in T(E)[t]/(t^{n+1}).$$

Note that $P_{\pi_1\pi_2}(t) = P_{\pi_1}(t)P_{\pi_2}(t)$ if ${\bold h}(\pi_1) = {\bold t}(\pi_2)$.
Write $$P_{\pi}(t) = \sum_{k=0}^{l(\pi)} (-1)^k e(\pi,k)t^k .$$
Set $e(\pi,k) = 0$ if $k> l(\pi).$
For $v \in \cup_{i=1}^n V_i$, set $P_v(t) = P_{\pi_v}(t)$ and $e(v,k) = e(\pi_v,k).$ Also, 
set $P_*(t) = 1$ and $e(*,k)= 0$ if $k > 0.$
\bigskip
{\bf Definition 2.1} Let $R$ be the ideal in $T(E)$ generated by 
$$\{ e(\pi_1,k) - e(\pi_2,k) \ | \ {\bold t}(\pi_1)={\bold t}(\pi_2), {\bold h}(\pi_1)={\bold h}(\pi_2), \ 1 \leq k \leq l(\pi_1) \}.$$

Note that this implies
$$P_{\pi_1}(t) \equiv P_{\pi_2}(t) \ \ mod \ R[t].$$
Now assume $v > u$, so  there is a path $\pi$ from $v$ to $u$.
Then $${\bold t}(\pi\pi_u) = v = {\bold t}(\pi_v) \ ,$$ $${\bold h}(\pi\pi_u) = * = {\bold h}(\pi_v)$$  
and so $$P_{\pi\pi_u}(t) \equiv P_{\pi_v}(t) \ \ mod \ R[t].$$
But  $$P_{\pi\pi_u}(t) = P_{\pi}(t)P_{\pi_u}(t) = P_{\pi}(t)P_u(t) \ ,$$
$$P_{\pi_v}(t) = P_v(t)$$  so   $$P_{\pi}(t) \equiv P_v(t)P_u(t)^{-1} \ mod \ R[t].$$ 

Noting that $P_{\pi}(t)$ is a polynomial of degree $l(\pi) = |v| - |u|$ and writing $(1 - a)^{-1} = 1 + a + a^2  ... \ $ for a nilpotent element $a$, we obtain
$$P_{\pi}(t) \equiv \sum_{r \geq 0,  i_0 \geq 0, \atop i_1,...,i_r \geq 1} (-1)^{i_0+...+i_r+r} e(v,i_0)e(u,i_1)...e(u,i_r)t^{i_0+...+i_r}$$
$$\equiv \sum_{j=0}^{l(\pi)} (\sum_{{{r \geq 0, i_0 \geq 0,} \atop {i_1,...,i_r \geq 1,}} \atop {i_0+...+i_r=j}} (-1)^{j+r} e(v,i_0)e(u,i_1)...e(u,i_r))t^j \ \ mod \ R[t].$$

Let$$H(v,u,j) =
\sum_{{{r \geq 0, i_0 \geq 0,} \atop {i_1,...,i_r \geq 1,}} \atop {\ i_0+...+i_r=j}} (-1)^{j+r} e(v,i_0)e(u,i_1)...e(u,i_r)$$
so that $$P_{\pi}(t) \equiv \sum_{j=0}^{l(\pi)} H(v,u,j)t^j \ \ mod \ R[t].$$

It follows that, modulo $R[t]$, 
$$P_v(t) \equiv  P_{\pi}(t)P_u(t)$$
$$ \equiv (\sum_{j=0}^{l(\pi)} H(v,u,j)t^j)(\sum_{i_{r+1}=0}^{|u|} (-1)^{i_{r+1}}e(u,i_{r+1})t^{i_{r+1}})$$

$$\equiv \sum_{{{r \geq 0, i_0,i_{r+1} \geq 0,} \atop {i_1,...,i_r \geq 1,}} \atop {{i_0+...+i_r \leq |v|-|u|,} \atop i_0+...+i_{r+1} \leq |v|}} (-1)^{i_0+...+i_{r+1}+r} e(v,i_0)e(u,i_1)...e(u,i_{r+1})t^{i_0+...+i_{r+1}} \ \mod \ R[t].$$

Setting $k = |v| - |u|$ and comparing coefficients of $t^{k+l}$ gives 
$$e(v,k+l) \equiv e(v,k)e(u,l) \ + $$
$$\sum_{{{r \geq 0, i_0,i_{r+1} \geq 0,} \atop  {i_1,...,i_r \geq 1,}} \atop {{i_0 < k,  i_0+...+i_r \leq k} \atop i_0+...+i_{r+1} = k + l}} (-1)^r e(v,i_0)e(u,i_1)...e(u,i_{r+1}) \ \ mod \ R . \tag{2.1}$$


Writing $$E(v,u,k,l) =
\sum_{{{r \geq 0, i_0,i_{r+1} \geq 0,} \atop  {i_1,...,i_r \geq 1,}} \atop {{i_0 < k,  i_0+...+i_r \leq k} \atop i_0+...+i_{r+1} = k + l}} (-1)^r e(v,i_0)e(u,i_1)...e(u,i_{r+1}) \ \ $$
we obtain  $$e(v,k+l) - e(v,k)e(u,l) \equiv E(v,u,k,l) \ \ mod \ R$$
when $v > u$, $|v| - |u| = k.$

We also have

\smallskip
\proclaim{Lemma 2.1} $E(v,u,1,l) = - e(u,1)e(u,l) + e(u,l+1)$.
\endproclaim
\smallskip
{\bf Proof:} 
Setting $k=1$ in the sum defining $E(v,u,k,l)$ 
gives $i_0 = 0, r = 0$ (hence $i_1 = l + 1$) or $i_0 = 0, r = 1$ 
(hence $i_1 = l, i_2 = l$). These two choices  give the two terms on the right-hand side.

\bigskip
Note that $|e(v,k)| = k|v| - k(k-1)/2$, and so, if $k = |v| - |u|$, 
$$|e(v,k)e(u,l)| = |e(v,k)| + |e(u,l)|$$
$$= k|v| - k(k-1)/2 + l|u| - l(l-1)/2$$
$$= (k + l)|v| - (k + l)(k + l - 1)/2$$
$$= |e(v, k + l)|.$$ 

We also have  
\bigskip
\proclaim{Lemma 2.2} If $|v| - |u| = k$,  $|E(v,u,k,l)| < |e(v,k+l)|$.
\endproclaim

{\bf Proof:} It is sufficient to show that each summand in the 
expression for $E(v,u,k,l)$ belongs to $T(E)_{|e(v,k+l)|-1}.$
Thus we must show  
$$ i_0|v| + (k + l - i_0)|u| - i_0(i_0 - 1)/2 - ... -i_{r + 1}(i_{r + 1} - 1)/2 < (k+l)|v| - (k+l)(k+l-1)/2,$$ 
where  
$r \ge 0, i_0 ,i_{r + 1} \geq 0$, $i_1, ..., i_r \geq 1$, $i_0 + ... + i_{r} \leq k$, $i_0 < k$  and 
$i_0 + ... + i_{r + 1} = k + l$.
Since $|u| = |v|-k$ this is equivalent to
$$-(k+l-i_0)k - (i_0^2 + ... + i_{r+1}^2 - k - l)/2 < -(k+l)(k+l-1)/2$$ 
which may be simplified to
$(k-i_0)^2 + ... + i_{r+1}^2 > l^2.$  This holds since $i_{r+1} \le l$ and $i_0 < k.$ 
\bigskip  

\proclaim{Lemma 2.3} a) If $f \in E$ then $f - e({\bold t}(f),1) + e({\bold h}(f),1) \in R$. 

b) If $v > u$, $|u| = |v| - 1$ then $e(v,1)e(u,k) - e(v,k+1) + e(u,k+1) - e(u,1)e(u,k) \in R$.
\endproclaim

{\bf Proof:} a) $P_{{\bold t}(f)}(t) \equiv P_{f\pi_{{\bold h}(f)}}(t) = (1 - tf)P_{{\bold h}(f)}(t) \ \ mod \ R[t],$ 
so $e({\bold t}(f),1) \equiv  f + e({\bold h}(f),1)  \ \ mod \ R$.

\smallskip

b) Since  $v > u$ and  $|u| = |v| - 1$ there is $f \in E$ with ${\bold t}(f) = v$, ${\bold h}(f) = u$. Then 
$$P_v(t) \equiv P_{f\pi_{{\bold h}(f)}}(t) = (1 - tf)P_u(t)$$
$$ \equiv (1 + te(v,1) - te(u,1))P_u(t) \ \ mod \ R[t].$$
Thus $$e(v,k+1) \equiv e(u,k+1) - e(v,1)e(u,k) + e(u,1)e(u,k) \ \ mod \ R .$$

\bigskip

\proclaim{Lemma 2.4} $e(v,1)E(u,w,k,l)- e(u,1)E(u,w,k,l) \equiv E(v,w,k+1,l)-E(u,w,k+1,l) \ \ $ $mod \ R .$
\endproclaim

{\bf Proof:} By the definition, the left-hand side (LHS) is 
$$\sum (-1)^r  e(v,1)e(u,i_0)e(w,i_1)...e(w,i_{r+1}) 
 - \sum (-1)^r e(u,1)e(u,i_0)e(w,i_1)...e(w,i_{r+1}), $$ 
where the sums are over 
$r \geq 0$, $i_0,i_{r+1} \geq 0$, $i_1,...,i_r \geq 1$, $i_0 < k$,  $i_0+...+i_r \leq k$ and $i_0+...+i_{r+1} = k + l.$

By Lemma 2.3,  $e(v,1)e(u,i_0) - e(u,1)e(u,i_0) \equiv   e(v,i_0+1) - e(u,i_0+1) \ \ mod \  R$.
Thus  the LHS is congruent to 
$$\sum (-1)^r  e(v,i_0)e(w,i_1)...e(w,i_{r+1}) 
 - \sum (-1)^r e(u,i_0)e(w,i_1)...e(w,i_{r+1}),$$ 
where the summation  is over all 
$r \geq 0$, $i_0 \geq 1$, $i_{r+1} \geq 0$, $i_1,...,i_r \geq 1$, $i_0 < k+1$,  $i_0+...+i_r \leq k+1$ and $i_0+...+i_{r+1} =  k + l+1.$

But this is equal to the same expression where the sum is over all 
$r \geq 0$, $i_0 \geq 0$, $i_{r+1} \geq 0$, $i_1,...,i_r \geq 1$, $i_0 < k+1$,  $i_0+...+i_r \le k+1$ and $i_0+...+i_{r+1} =  k + l+1$ (since $e(v,0) = e(u,0) = 1.$) This is the right hand side of the asserted congruence. 
 
\bigskip
\proclaim{Lemma 2.5} $R$ is generated by all $e(\pi_1,k) - e(\pi_2,k)$, ${\bold t}(\pi_1) = {\bold t}(\pi_2)$, ${\bold h}(\pi_1) = {\bold h}(\pi_2) = *$.
\endproclaim

{\bf Proof:} Let $S$ be the ideal generated by all such elements. Thus, for such $\pi_1$, $\pi_2$, 
$$ P_{\pi_1}(t) \equiv P_{\pi_2}(t) \ \ mod \ S[t].$$
Let $P_{\pi_3}(t), P_{\pi_4}(t) \in T(E)[t]/(t^{n+1})$  satisfying ${\bold t}(\pi_3) = {\bold t}(\pi_4)$, ${\bold h}(\pi_3) = {\bold h}(\pi_4) = w$.  Then 
${\bold t}(\pi_3\pi_w) = {\bold t}(\pi_3) = {\bold t}(\pi_4) = {\bold t}(\pi_4\pi_w)$, ${\bold h}(\pi_3\pi_w) = {\bold h}(\pi_w) = * = {\bold h}(\pi_4\pi_w)$ and 
$$P_{\pi_3}(t)P_{\pi_w}(t) = P_{\pi_3\pi_w}(t) \equiv P_{\pi_4\pi_w}(t) = P_{\pi_4}(t)P_{\pi_w}(t) \ \ mod \ S[t].$$
Then since $P_{\pi_w}(t)$ is invertible,  $P_{\pi_3}(t) \equiv P_{\pi_4}(t) \ \ mod \ S[t]$.   Since the coefficients of all $P_{\pi_3}(t) - P_{\pi_4}(t)$ generate $R$,  
we have  the result.

\bigskip
\proclaim{Lemma 2.6} Let $S_1 = \{f - e({\bold t}(f),1) + e({\bold h}(f),1)| f \in E\}$ and 
$$S_2 = \{- e(v,1)e(u,k) - e(v,k+1) + e(u,k+1) + e(u,1)e(u,k)|u,v \in V, v>u, |u| = |v| - 1 > 0\}.$$ 
Then $S_1 \cup S_2$ generates $R$.
\endproclaim

{\bf Proof:} By Lemma 2.3, $S_1 \cup S_2 \subseteq R.$  Let $v \in V, |v| > 0$ and let
$\pi = \{e_1,e_2,...,e_{|v|}\}$ be a path from $v$ to $*$.  By lemma 2.5, it is sufficient to show that
$P_{\pi}(t) \equiv P_v(t) \ mod \ (S_1 \cup S_2)[t].$  We proceed by induction on $|v|.$  If $|v| = 1$, 
then $\pi = \{e_1\}$ where ${\bold t}(e_1) = v, {\bold h}(e_1) = *.$  Then
$P_{\pi}(t) = 1 - te_1 \equiv 1 + e(v,1) - e(*,1) \equiv 1 - te_v \equiv P_v(t) \ mod \ 
S_1[t].$

Now assume $|v| > 1$ and write $\pi = e_1\pi'$ with ${\bold t}(e_1) = v, {\bold h}(e_1) = u.$  Then $|u| < |v|$ and so by induction we have
$P_u(t) \equiv P_{\pi'}(t) \ mod \ (S_1 \cup S_2)[t].$  
Since $$P_v(t) \equiv (1 - te(v,1) + te(u,1))P_u(t)\ mod \ S_2[t]$$ and $$
1 - te(v,1) + te(u,1) \equiv 1 - te_1 \ mod S_1[t]$$ we have
$$P_v(t) \equiv (1 - te_1)P_u(t) \equiv (1 - te_1)P_{\pi'}(t) \equiv P_{\pi}(t) \ mod \ (S_1 \cup S_2)[t],$$
as required.
\bigskip

\head \bf 3. The algebra $A(\Gamma)$
\endhead
\bigskip
Let $$A(\Gamma) = T(E)/R$$
and
$$A(\Gamma)_i = (T(E)_i + R)/R  .$$

This gives $A(\Gamma)$ the structure of a filtered algebra. 
Let $\tilde{}$ denote the canonical homomorphism $\tilde{}: T(E) \rightarrow  T(E)/R = A(\Gamma)$ 
and write $\tilde{e}(v,k$), $\tilde{E}(v,u,k,l)$, etc., for the images of the elements 
$e(v,k)$, $E(v,u,k,l)$, etc. If $a \in A(\Gamma)_i, a \not\in A(\Gamma)_{i-1},$ write $|a| = i.$

Note that 
$$\{ \tilde{e}(v,k) \ | \ v \in \cup_{i=1}^n V_i, \  \  k \leq |v| \}$$
generates $A(\Gamma)$,  since by Lemma 2.3a), $\tilde{f} = \tilde{e}({\bold t}(f),1) - \tilde{e}({\bold h}(f),1)$ for 
any $f \in E$.   

We now develop some notation for products of the $\tilde{e}(v,k)$. 

We say that a pair $(v,k)$, $v \in V$, $0 \leq k \leq |v|$, can be {\it composed} with the pair $(u,l)$, 
$u \in V$, $0 \leq l \leq |u|$, if $v > u$ and $|u| = |v| - k$. If $(v,k)$ can be composed with $(u,l)$ 
we write $ (v,k) \models (u,l)$.

Let ${\bold B_1}(\Gamma)$   be the set of all sequences
$$ {\bold b} = ((b_1,m_1),(b_2,m_2),...,(b_k,m_k))$$  where 
$k \geq 0$, $b_1, b_2, ..., b_k \in V$, $0 \leq m_i \leq |b_i|$ for $1 \leq i \leq k$.
Let $\emptyset$ denote the empty sequence. Define $|{\bold b}|,$ the {\it level} of ${\bold b}$, by
$$|{\bold b}| = \sum_{i=1}^k \{m_i|b_i| - {{m_i(m_i-1)} \over 2}\}.$$

If $1 \leq s \leq k$ write
${\bold b}^s = ((b_s,m_s),...,(b_k,m_k))$. Write  ${\bold b}^{k+1} = \emptyset.$

If $ {\bold b} = ((b_1,m_1),(b_2,m_2),...,(b_k,m_k))$ and ${\bold c} = ((c_1,n_1),(c_2,n_2),...,(c_s,n_s)) \in 
{\bold B_1}(\Gamma)$ define 
$$ {\bold b} \circ {\bold c} = ((b_1,m_1),(b_2,m_2),...,(b_k,m_k),(c_1,n_1),(c_2,n_2),...,(c_s,n_s)).$$

Let $${\bold B}(\Gamma) = 
 \{  {\bold b} = ((b_1,m_1),(b_2,m_2),...,(b_k,m_k)) \in {\bold B_1}(\Gamma) \ | \ \atop 
        (b_i,m_i) \not\models  (b_{i+1},m_{i+1}) ,   \  \  1  \leq i < k \}.$$

For  $${\bold b} = ((b_1,m_1),(b_2,m_2),...,(b_k,m_k)) \in {\bold B_1}(\Gamma)$$ 
set 
$$\tilde{e}({\bold b}) = \tilde{e}(b_1,m_1)...\tilde{e}(b_k,m_k).$$
Then $|\tilde{e}({\bold b})| = |{\bold b}|$ 
and $\tilde{e}({\bold b} \circ {\bold c}) = \tilde{e}({\bold b})\tilde {e} ({\bold c}).$
Clearly
$\{ \tilde{e}({\bold b}) \ | \ {\bold b} \in {\bold B_1}(\Gamma) \}$ spans $A(\Gamma)$.
The following lemma is immediate from (2.1) and the definition of $E(v,u,k,l).$

\bigskip
\proclaim{Lemma 3.1} If $(v,k) \models (u,l)$ then $\tilde{e}(v,k)\tilde{e}(u,l) = \tilde{e}(v,k+l)
- \tilde{E}(v,u,k,l)$.
\endproclaim

For ${\bold b} = ((b_1,m_1),(b_2,m_2),...,(b_k,m_k)) \in {\bold B}(\Gamma)$ 
define $$z(v,k,{\bold b}) = min \{\{ j \ | \ (v,k+m_1+...+m_{j-1}) \not\models (b_j,m_j) \}\cup \{k+1\}\}.$$

\bigskip
\proclaim{Lemma 3.2} Let ${\bold b} \in {\bold B}(\Gamma)$. Then $\tilde{e}(v,k)\tilde{e}({\bold b}) = $
$$\tilde{e}((v,k+m_1+...+m_{z(v,k,{\bold b})-1)}) \circ {\bold b}^{z(v,k,{\bold b})})$$
$$- \sum_{j=1}^{z(v,k,{\bold b})-1} \tilde{E}(v,b_j,k+m_1+...+m_{j-1},m_j)\tilde{e}({\bold b}^{j+1}).$$
\endproclaim

\smallskip
{\bf Proof:} For $1 \le i \le z(v,k,{\bold b})-1$ , we have $(v,k+m_1+...+m_{i-1})\models (b_i,m_i)$ and so, by Lemma 3.1,
$\tilde {e}(v,k+m_1+...+m_{i-1})\tilde{e}({\bold b}^i) = \tilde{e}((v,k+m_1+...+m_{i-1}+m_i)\circ{\bold b}^{i+1}) - E(v,b_i,k+...+m_{i-1},m_i).$  It follows, by induction on $i$, that for $1 \leq i \leq   z(v,k,{\bold b})-1$ 

$$\tilde{e}(v,k)\tilde{e}({\bold b}) = $$
$$\tilde{e}((v,k+m_1+...+m_{i-1}+m_i) \circ {\bold b}^{i+1})$$
$$- \sum_{j=1}^{i} \tilde{E}(v,b_j,k+m_1+...+m_{j-1},m_j)\tilde{e}({\bold b}^{j+1}).$$
Taking $i = z(v,k,{\bold b}) - 1$ gives the lemma.

\bigskip
\proclaim{Corollary 3.3} $S = \{ \tilde{e}({\bold b}) \ | \ {\bold b} \in {\bold B}(\Gamma) \}$ spans $A(\Gamma).$
\endproclaim

{\bf Proof:} Since $1 \in S$, it is sufficient to show  that $span S$ is invariant under multiplication 
by $\tilde{e}(v,k)$, hence sufficient to show that $\tilde{e}(v,k)\tilde{e}({\bold b}) \in span S$, if 
${\bold b} \in {\bold B}(\Gamma)$.  We prove this by induction on $|\tilde{e}(v,k)\tilde{e}({\bold b})|$.
If $|\tilde{e}(v,k)\tilde{e}({\bold b}) | = 0$, then $\tilde{e}(v,k)\tilde{e}({\bold b}) \in F1 \subseteq \ span \  S$.  
Now assume that $\tilde{e}(w,l)\tilde{e}({\bold c}) \in \ span \ S $ whenever 
${\bold c} \in {\bold B}(\Gamma)$ and $|\tilde{e}(w,l)\tilde{e}({\bold c})| < |\tilde{e}(v,k)\tilde{e}({\bold b})|.$  In 
view of Lemma 2.2, this implies $\tilde{E}(v,b_1,k+m_1+...+m_{j-1},m_j)\tilde{e}({\bold b}^{j+1}) \in S$ for
$1 \le j \le z(v,k,{\bold b}) -1$.  But by Lemma 3.2,

$$\tilde{e}(v,k)\tilde{e}({\bold b}) = $$
$$\tilde{e}((v,k+m_1+...+m_{z(v,k,{\bold b})-1}) \circ {\bold b}^{z(v,k,{\bold b})})$$
$$- \sum_{j=1}^{z(v,k,{\bold b})-1} \tilde{E}(v,b_j,k+m_1+...+m_{j-1},m_j)\tilde{e}({\bold b}^{j+1}).$$
Since $(v,k+m_1+...+m_{z(v,k,{\bold b})-1)}\circ{\bold b}^{z(v,k,{\bold b})} \in {\bold B}(\Gamma),$
every summand on the right-hand side of this expression belongs to $S$, so 
$\tilde {e}(v,k)\tilde {e}({\bold b}) \in S$, as required.
\bigskip

\head \bf 4.  Independence Theorem 
\endhead
\bigskip
Define $B$ to be the vector space over $F$ with basis ${\bold B}(\Gamma)$.  Let 
$$B_i = span \{ {\bold b} \ | \ |{\bold b}| \leq i \}.$$ 

We will define a linear transformation
$$\mu: T(E) \otimes B \rightarrow B$$
giving $B$ the structure of a $T(E)$ module.

First define $$\mu_1: E \times {\bold B}(\Gamma) \rightarrow B$$ by
$$\mu_1: (f,\emptyset) \mapsto e({\bold t}(f),1)$$
for $f \in E_1$ (where $\emptyset$ denotes the empty sequence) and
$$\mu_1: E_i \times B_h \rightarrow 0$$
if $i + h >1.$  Then, as $T(E)$ is the free algebra on the set $E$, 
$\mu_1$ extends to a linear transformation, again denoted $\mu_1$ from 
$T(E) \otimes B$ to $B$ giving $B$ the structure of a $T(E)$- module.

Now assume $s > 1$ and that we have defined maps 
$$\mu_j:E \times {\bold B}(\Gamma) \rightarrow B$$
for $1 \le j \le s-1$ such that
$$\mu_j: E_i \times {\bold B}(\Gamma)_h \rightarrow B_{i+h}$$
for all $1 \le i \le n$ and $h \ge 0$.  Assume also that
$$\mu_j: E_i \times {\bold B}(\Gamma)_h \rightarrow 0 $$
whenever $i+h \ge j$ and that 
$$\mu_{j'}|_{E_i \times {\bold B}(\Gamma)_h} = \mu_{j''}|_{E_i \times {\bold B}(\Gamma)_h}$$
whenever $j' \ge j'' \ge i+h.$ As in the case of $\mu_1$, each $\mu_j$ extends to a linear transformation, again denoted 
$\mu_j$ from $T(E) \otimes B$ to $B$ giving $B$ the structure of a $T(E)$-module.

We now define 
$$\mu_s: E \times {\bold B}(\Gamma) \rightarrow B$$ by
$$\mu_s|_{E_i \times {\bold B}(\Gamma)_h} = \mu_{s-1}|_{E_i \times {\bold B}(\Gamma)_h}$$
for $i+h \ne s$, and  

$$\mu_s: (f,{\bold b}) \mapsto   (v,1+m_1+...+m_{z(v,1,{\bold b})-1}) \circ {\bold b}^{z(v,1,{\bold b})}$$

$$- \sum_{j=1}^{z(v,1,{\bold b})-1} E(v,b_j,1+m_1+...+m_{j-1},m_j){\bold b}^{j+1} \tag{4.1}$$
$$- (u,1+m_1+...+m_{z(u,1,{\bold b})-1}) \circ {\bold b}^{z(u,1,{\bold b})}$$
$$+ \sum_{j=1}^{z(u,1,{\bold b})-1} E(u,b_j,k+m_1+...+m_{j-1},m_j){\bold b}^{j+1}$$
for ${\bold b} = ((b_1,m_1),...,(b_k,m_k)), {\bold t}(f) = v, {\bold h}(f) = u, |f| + |{\bold b}| = s.$ (Note that, as
$e(*,j) = 0$ for $j > 0$, the last two summands vanish when $|f| = 1$.) 

Thus we have inductively defined $\mu_j$ for all $j$.  Define
$$\mu: T(E) \otimes B \rightarrow B$$
by
$$\mu|_{T(E)_i \otimes B_h} = \mu_{i+h}|_{T(e)_i \otimes B_h}.$$

To simplify the notation we write $f{\bold b}$ for $\mu(f \otimes {\bold b}).$

\bigskip
\proclaim{Lemma 4.1} (a) $e(v,k){\bold b} = (v,k+m_1+...+m_{z(v,k,{\bold b})-1}) \circ {\bold b}^{z(v,k,{\bold b})}$
$$- \sum_{j=1}^{z(v,k,{\bold b})-1}E(v,b_j,k+m_1+...+m_{j-1},m_j){\bold b}^{j+1}.$$

(b) $RB = (0)$
\endproclaim

{\bf Proof:} Note that (a) is clear if $k = 0$ (for $e(v,0) = 1$ and $z(v,0,{\bold b}) = 1$).  
Also (4.1) shows that (a) holds when $k = 1$ (as $e(v,1) = e_{v^{(0)}} + e_{v^{(1)}} + ... + e_{v^{(|v|-1)}}$ and so the expresion given by (4.1) for $e(v,1){\bold b}$ is a telescoping series). Thus (a) holds whenever $|e(v,k)| + |{\bold b}| \le 1.$

Let $R_i = R \cap A(\Gamma)_i$ and note that $R_0 = (0)$ and, by Lemma 2.6, $R_1$ is spanned by
$\{f - e({\bold t}(f),1)|f \in E_1\}.$  Hence $R_1B_0$ is spanned by
$\{f\emptyset - e({\bold t}(f),1)\emptyset|f \in E_1\} = \{0\},$ so $R_iB_h = (0)$ whenever $i + h \le 1.$
 
Now assume $s > 1$ and that $$e(w,l){\bold c} = (w,k+n_1+...+n_{z(w,l,{\bold c})-1}) \circ {\bold c}^{z(w,l,{\bold c})}$$
$$- \sum_{j=1}^{z(w,l,{\bold c})-1}E(w,c_j,l+n_1+...+n_{j-1},n_j){\bold c}^{j+1}$$
whenever ${\bold c} = ((c_1,n_1)...(c_p,c_p)) \in {\bold B}(\Gamma)$. Assume also that
$|e(w,l)| + |{\bold c}| < s$ and $R_iB_h = (0)$ whenever $i + h < s$.  We will show that if $|e(v,k)| + |{\bold b}| = s$, then (a) holds and that if $i + h = s$ then $R_iB_h = (0)$, thus proving the lemma by induction.

By the expression obtained for $e(v,1){\bold b}$, we have 
$$(f - e({\bold t}(f),1) + e({\bold h}(f),1)){\bold b} = 0$$
for all $f \in E, {\bold b} \in {\bold B}(\Gamma).$ Thus, by Lemma 2.6, to show $R_iB_h = (0)$ whenever 
$i + h = s$ it is sufficient to show that
$$(- e(v,1)e(u,k) - e(v,k+1) + e(u,k+1) + e(u,1)e(u,k)){\bold b} = (0)$$  
whenever $|e(v,k)| + |{\bold b}| = s.$

\smallskip

Now assume that $|e(v,k)| + |{\bold b}| = s$.  Since $|e(u,k-1)| < |e(v,k)|$ and $|e(u,k)| < |e(v,k)|$,  the 
induction assumption gives values for $e(u,k-1){\bold b}$ and $e(u,k){\bold b}$.  Thus 

$$e(v,1)e(u,k-1){\bold b} + e(u,k){\bold b} - e(u,1)e(u,k-1){\bold b}$$
$$ = e(v,1)(u,k-1+m_1+...+m_{z(u,k-1,{\bold b})-1}) \circ {\bold b}^{z(u,k-1,{\bold b})}$$
$$- \sum_{j=1}^{z(u,k-1,{\bold b})-1}e(v,1)E(u,b_j,k-1+m_1+...+m_{j-1},m_j){\bold b}^{j+1}$$

$$+ (u,k+m_1+...+m_{z(u,k-1,{\bold b})-1}) \circ {\bold b}^{z(u,k,{\bold b})}$$
$$- \sum_{j=1}^{z(u,k-1,{\bold b})-1}e(v,1)E(u,b_j,k+m_1+...+m_{j-1},m_j){\bold b}^{j+1}$$

$$- e(u,1)(u,k-1+m_1+...+m_{z(u,k-1,{\bold b})-1}) \circ {\bold b}^{z(u,k-1,{\bold b})}$$
$$+ \sum_{j=1}^{z(u,k-1,{\bold b})-1}e(u,1)E(u,b_j,k-1+m_1+...+m_{j-1},m_j){\bold b}^{j+1}$$

\noindent  By the definition of the actions of $e(v,1)$ and $e(u,1)$ this becomes 

$$(v,k+m_1+...+m_{z(v,k,{\bold b})-1}) \circ {\bold b}^{z(u,k-1,{\bold b})}$$
$$-E(v,u,1,k-1+m_1+...+m_{z(v,k,{\bold b})-1}){\bold b}^{z(u,k-1,{\bold b})}$$
$$- \sum_{j=z(u,k-1,{\bold b})-1}^{z(v,k,{\bold b})-1}E(v,b_j,k+m_1+...+m_{j-1},m_j){\bold b}^{j+1}$$
$$- \sum_{j=1}^{z(u,k-1,{\bold b})-1}e(v,1)E(u,b_j,k-1+m_1+...+m_{j-1},m_j){\bold b}^{j+1}$$

$$+ (u,k+m_1+...+m_{z(u,k,{\bold b})-1}) \circ {\bold b}^{z(u,k,{\bold b})}$$
$$-  \sum_{j=1}^{z(u,k,{\bold b})-1}E(u,b_j,k+m_1+...+m_j){\bold b}^{j-1}$$

$$- (u,1) \circ (u,k-1+m_1+...+m_{z(u,k-1,{\bold b})-1}) \circ {\bold b}^{z(u,k-1,{\bold b})}$$
$$+ \sum_{j=1}^{z(u,k-1,{\bold b})-1} e(u,1)E(u,b_j,k-1+m_1+...+m_{j-1},m_j) {\bold b}^{j+1}.$$

\noindent Set $$G = - e(v,1)E(u,b_j,k-1+m_1+...+m_{j-1},m_j)$$
$$+ e(u,1)E(u,b_j,k-1+m_1+...+m_{j-1},m_j)$$
$$+ E(v,b_j,k+m_1+...+m_{j-1},m_j)$$
$$- E(u,b_j,k+m_1+...+m_j).$$
\noindent By Lemma 2.6, $G \in R$  and by Lemma 2.2 $|G{\bold b}^{j+1}| < s$. Hence by the induction assumption, 
$G{\bold b}^{j+1} = 0$.  Thus we may replace 
$$ e(v,1)E(u,b_j,k-1+m_1+...+m_{j-1},m_j){\bold b}^{j+1}$$
$$- e(u,1)E(u,b_j,k-1+m_1+...+m_{j-1},m_j){\bold b}^{j+1}$$
by
$$ E(v,b_j,k+m_1+...+m_{j-1},m_j){\bold b}^{j+1}$$
$$- E(u,b_j,k+m_1+...+m_j){\bold b}^{j+1}.$$

Thus our expression becomes  
$$(v,k+m_1+...+m_{z(v,k,{\bold b})-1}) \circ {\bold b}^{z(u,k-1,{\bold b})}$$
  
$$-E(v,u,1,k-1+m_1+...+m_{z(v,k,{\bold b})-1}){\bold b}^{z(u,k-1,{\bold b})}$$

$$- \sum_{j=1}^{z(v,k,{\bold b})-1} E(v,b_j,k+m_1+...+m_{j-1},m_j){\bold b}^{j+1}$$

$$+ (u,k+m_1+...+m_{z(u,k,{\bold b})-1}) \circ {\bold b}^{z(u,k,{\bold b})}$$
$$- (u,1) \circ (u,k-1+m_1+...+m_{z(u,k-1,{\bold b})-1}) \circ {\bold b}^{z(u,k-1,{\bold b})}$$

\noindent But by Lemma 2.1 $$E(v,u,1,k-1+m_1+...+m_{z(v,k,{\bold b})-1}) $$ 
$$= - e(u,1)e(u,k-1+m_1+...+m_{z(v,k,{\bold b})-1}) + e(u,k+m_1+...+m_{z(v,k,{\bold b})-1}).$$
\noindent Since 
$$|E(v,u,1,k-1+m_1+...+m_{z(v,k,{\bold b})-1})| + |{\bold b}^{z(u,k-1,{\bold b})}| < s,$$
by Lemma 2.2, and (as $|u| < |v|$)
$$|e(u,1)e(u,k-1+m_1+...+m_{z(v,k,{\bold b})-1})| + |{\bold b}^{z(u,k-1,{\bold b})}| < s,$$
and
$$|e(u,k+m_1+...+m_{z(v,k,{\bold b})-1})| + |{\bold b}^{z(u,k-1,{\bold b})}| < s,$$
the induction assumption allows us to replace 
$$E(v,u,1,k-1+m_1+...+m_{z(v,k,{\bold b})-1}) {\bold b}^{z(u,k-1,{\bold b})}$$
\noindent by
$$- (u,1) \circ (u,k-1+m_1+...+m_{z(u,k-1,{\bold b})-1}) \circ {\bold b}^{z(u,k-1,{\bold b})}$$

$$+ (u,k+m_1+...+m_{z(u,k,{\bold b})-1}) \circ {\bold b}^{z(u,k,{\bold b})}.$$

Making this substitution gives 
$$e(v,1)e(u,k-1){\bold b} + e(u,k){\bold b} - e(u,1)e(u,k-1){\bold b}$$

$$= (v,k+m_1+...+m_{z(v,k,{\bold b})-1}) \circ {\bold b}^{z(v,k,{\bold b})} \tag{4.2}$$

$$- \sum_{j=1}^{z(v,k,{\bold b})-1} E(v,b_j,k+m_1+...+m_{j-1},m_j){\bold b}^{j+1}.$$

Since $P_v(t) = (1 + te(v,1) - te(v^{(1)},1))P_{v^{(1)}}(t)$ we have
$$e(v,k) = e(v,1)e(v^{(1)},k-1) + e(v^{(1)},k) - e(v^{(1)},1)e(v^{(1)},k-1).$$
Thus setting $u = v^{(1)}$ in (4.2) gives (a). Since the right-hand side of (4.2) is independent of $u$ 
we also obtain  
$$e(v,k){\bold b} = e(v,1)e(u,k-1){\bold b} + e(u,k){\bold b} - e(u,1)e(u,k-1){\bold b}$$
for all $u$ with $v > u, |u| = |v| - 1.$  As noted above, this completes the proof of (b).

\bigskip
\proclaim{Corollary 4.2} There is an action of $A(\Gamma)$ on $B$ satisfying 
$$\tilde{e}(v,k){\bold b} = (v,k+m_1+...+m_{z(v,k,{\bold b})-1}) \circ {\bold b}^{z(v,k,{\bold b})}$$

$$- \sum_{j=1}^{z(v,k,{\bold b})-1} \tilde{E}(v,b_j,k+m_1+...+m_{j-1},m_j){\bold b}^{j+1}.$$
\endproclaim
\bigskip
Consequently, for ${\bold b} \in {\bold B}(\Gamma)$, $e({\bold b})1 = {\bold b}$.
Since  ${\bold B}(\Gamma) \subset B$  is linearly independent,  
$\{ \tilde{e}({\bold b}) \ | \  {\bold b} \in {\bold B}(\Gamma) \}$  is linearly independent.  
Therefore,  we have proved the following theorem:

\bigskip  

\proclaim{Theorem 4.3} Let $\Gamma = (V,E)$  be a layered graph, $V = \cup_{i=0}^n V_i, $ and $V_0 = \{*\}$ 
where $*$ is the unique minimal vertex of $\Gamma$. Then $\{ \tilde{e}({\bold b}) \ | \  {\bold b} \in {\bold B}(\Gamma) \}$ is a basis 
for $A(\Gamma)$.
\endproclaim
\medskip
Note that 
$$\{\tilde{e}({\bold b})| |\tilde {e}({\bold b})| \le i\}$$
is a basis for $A(\Gamma)_i$.  Therefore, writing
$\bar {e}({\bold b}) + A(\Gamma)_{i-1} \in gr \ A(\Gamma)$
where $|\tilde {e}({\bold b})| = i$ we have:
\medskip
\proclaim{Corollary 4.4} Let $\Gamma = (V,E)$  be a layered graph, $V = \cup_{i=0}^n V_i, $ and $V_0 = \{*\}$ 
where $*$ is the unique minimal vertex of $\Gamma$. Then $\{\bar{e}({\bold b})| {\bold b} \in {\bold B}(\Gamma)\}$
is a basis for $gr \ A(\Gamma).$
\endproclaim

Also, if $|\tilde {e} (v,k)| = i,$ we have
$$\tilde {e}(v,k) + A(\Gamma)_{i-1} = \tilde {e}_{v^{(0)}}\tilde {e}_{v^{(1)}}...
\tilde {e}_{v^{(k-1)}} + A(\Gamma)_{i-1}.$$
Write 
$$\hat {e}(v,k) = \hat {e}_{v^{(0)}}\hat {e}_{v^{(1)}}...
\hat {e}_{v^{(k-1)}}$$ and, for ${\bold b} = ((v_1,k_1),...,(v_s,k_s)) \in {\bold B}(\Gamma),$ write
$$\hat {e}({\bold b}) = \hat {e}(v_1,k_1)...\hat {e}(v_s,k_s).$$
Then, if $|\tilde {e}({\bold b}) | = i$, it follows that
$$\tilde {e}({\bold b}) + A(\Gamma)_{i-1} = \hat {e}({\bold b}) + A(\Gamma)_{i-1}$$
and so we have:
\medskip
\proclaim{Corollary 4.5} Let $\Gamma = (V,E)$  be a layered graph, $V = \cup_{i=0}^n V_i, $ and $V_0 = \{*\}$ 
where $*$ is the unique minimal vertex of $\Gamma$. Then $\{\hat{e}({\bold b})| {\bold b} \in {\bold B}(\Gamma)\}$
is a basis for $A(\Gamma).$
\endproclaim
Note that, if $\Gamma = \Gamma_n$, the basis for $Q_n = A(\Gamma_n) $ given by Corollary 4.5 is the basis constructed in \cite{GRW}.


\Refs \widestnumber\key{GKLLRT}

\ref\key GR1 \by I. Gelfand, V. Retakh \paper Noncommutative Vieta theorem 
and symmetric functions \book  Gelfand Mathematical Seminars 1993-95 \publ 
Birkhauser Boston \yr 1996 \pages 93-100 \endref

\ref\key GR2 \by I. Gelfand, V. Retakh \paper Quasideterminants I\jour 
Selecta Math. (N.S.) \vol 3 \yr 1997 \pages 517-546  \endref

\ref\key GGRSW \by I. Gelfand, S. Gelfand, V. Retakh, S. Serconek,
and R. Wilson\paper Hilbert series of quadratic algebras
associated with decompositions of noncommutative polynomials \jour
J. Algebra \vol 254 \yr 2002\pages 279--299  \endref

\ref\key GGRW \by I. Gelfand, S. Gelfand, V. Retakh, R.
Wilson \paper Quasideterminants \jour
Advances in Math. \vol 193  \yr 2005 \pages 56-141 \endref

\ref\key GRW \by I. Gelfand, V. Retakh, and R. Wilson \paper
Quadratic-linear algebras associated with decompositions of
noncommutative polynomials and Differential polynomials \jour
Selecta Math. (N.S.) \vol 7\yr 2001 \pages 493--523 \endref

\ref\key O  \by A. Odesskii \paper  Set-theoretical solutions to the Yang-Baxter relation
from factorization of matrix polynomials and $\theta$-functions \jour  
Mosc. Math. J. \vol 3 \issue 1 \yr 2003 \pages  97--103, 259 \endref

\ref\key Pi \by D. Piontkovski \paper  Algebras associated to pseudo-roots 
of noncommutative polynomials  are Koszul \paperinfo math.RA/0405375 \endref

\ref\key SW \by S. Serconek and R. L. Wilson \paper Quadratic algebras
associated with decompositions of noncommutative polynomials
are Koszul algebras
\jour J. Algebra \vol 278 \yr 2004 \pages 473-493 \endref


\ref\key V \by  A. Veselov \paper  Yang-Baxter maps and integrable dynamics \jour
Phys. Lett. A \vol 314 \issue 3 \yr 2003 \pages 214--221 \endref
\endRefs
\enddocument